\documentclass[preprint,3p,12pt]{elsarticle}

\usepackage{lineno,hyperref}
\modulolinenumbers[5]

\usepackage{hyperref}

\usepackage{epstopdf}

\usepackage{mathrsfs}
\usepackage{amssymb}
\usepackage{amsfonts}
\usepackage{latexsym}
\usepackage{amsmath}
\usepackage{xcolor}
\usepackage{wrapfig}
\usepackage{floatflt}

\usepackage{mathtools}
\usepackage{extarrows}

\usepackage{graphicx}
\def\lb{\label}

\newcommand{\er}[1]{\textrm{(\ref{#1})}}

\newtheorem{theorem}{\bf Theorem}[section]

\def\r{\rho}           
\def\s{\sigma}         
\def\S{\Sigma}

       \def\vp{\varphi}

    \def\N{{\mathbb N}}   

\def\lt{\biggl}                  \def\rt{\biggr}
               \def\wt{\widetilde}
\def\no{\noindent}

\let\ge\geqslant                 \let\le\leqslant

\def\iy{\infty}

\def\el2{\ell^{\,2}}             \def\1{1\!\!1}

\let\ge\geqslant
\let\le\leqslant

\newcommand{\ca}{\begin{cases}}
\newcommand{\ac}{\end{cases}}
\newcommand{\ma}{\begin{pmatrix}}
\newcommand{\am}{\end{pmatrix}}
\def\eq{\begin{equation}}
\def\qe{\end{equation}}
\def\[{\begin{equation}}
\def\]{\end{equation}}

\begin{document}

\begin{frontmatter}



\title{An entire function connected with the approximation of the golden ratio}

\date{\today}

\author
{Anton A. Kutsenko}
\address{Jacobs University, 28759 Bremen, Germany; email: akucenko@gmail.com}
\address{Saint-Petersburg State University,
Universitetskaya nab. 7/9, St. Petersburg, 199034, Russia}

\begin{abstract}
In 1987, R. B. Paris uses the analytic function 
\[\label{main}
 g(w)=\lim_{n\to\infty}(2\varphi)^n\biggl(\underbrace{\sqrt{1+\sqrt{1+...\sqrt{1+w}}}}_n-\varphi\biggr),\ \ \ \varphi=\frac{1+\sqrt{5}}2,
\]
to estimate the convergence of nested squares to the golden ratio. The function $g$ is non-entire and, perhaps, can not be expressed in terms of some standard known functions. We show that $f(z):=g^{-1}(z)$ is an entire function satisfying Poincare equality. While $f$ has zeros of various multiplicities, it can be expressed in terms of its simple zeros, forming fractal structures similar to Julia sets. 

\end{abstract}


\begin{keyword}
Golden ratio, nested squares, Poincare functions, Julia sets, Weierstrass-Hadamard expansion
\end{keyword}


\end{frontmatter}

{\section{Introduction}\lb{sec0}}

In his nice and short paper \cite{P}, R. B. Paris investigated the convergence of nested squares to the golden ratio. He found
$$
 \vp-\vp_n\sim\frac{2C}{(2\vp)^n},
$$
where $\vp_1=1$, $\vp_{n+1}=\sqrt{1+\vp_n}$. The value $C\approx1.098$ admits explicit representation $C=\vp\wt g(1/\vp)$, where
the analytic function $\wt g$ satisfies
$$
 \wt g(z)=2\vp\wt g(\vp-\sqrt{\vp^2-z}),\ \ |z|<\vp^2;\ \ \wt g(0)=0,\ \ \wt g'(0)=1,
$$
see \cite{F,W}. There is another form of this representation $2C=-g(0)$, see \er{main}. The functions $\wt g$, $g$ are analytic in some domain, but non-entire. We show that $f:=g^{-1}$ admits extension to an entire function. We explore some properties of $f$. Section \ref{sec1} contains some basic results, other sections are more experimental. 
Many of the results about Poincare functions have been known for a long time, see, e.g., \cite{EL,ES,DGV}. We try in a simple manner to show the connection between fractals, entire functions, and infinite expansions. Most of the general ideas about Poincare functions can be well understood by studying concrete examples. At the same time, examples may contain interesting expansions or representations not available in the general case. This note appear as a result of discussion I opened in \cite{mathforum}. In principle, I do not expect a great novelty in this note.

{\section{Basic properties of $g$ and $g^{-1}$}\lb{sec1}}
Let us note that \er{main} admits an infinite product expansion
\[\lb{main1}
 g(w)=(w-\vp)\cdot\frac{2\vp}{\vp+\sqrt{1+w}}\cdot\frac{2\vp}{\vp+\sqrt{1+\sqrt{1+w}}}\cdot...
\]
converging relatively fast as $\prod_n(1-(2C)(2\vp)^{-n})$.
Differentiating \er{main}, we obtain another infinite product expansion
\[\lb{difg}
 g'(w)=\frac{\vp}{\sqrt{1+w}}\cdot\frac{\vp}{\sqrt{1+\sqrt{1+w}}}\cdot\frac{\vp}{\sqrt{1+\sqrt{1+\sqrt{1+w}}}}....
\]
Thus $g'(\vp)=1$, since $w=\vp$ is a fixed point of the mapping $\sqrt{1+w}$. Along with $g(\vp)=0$, we can state that there is an analytic function $f(z)=g^{-1}(z)$ defined in some neighborhood of $z=0$. This function satisfies $f(0)=\vp$, $f'(0)=1$.

Definition \er{main} leads to the functional equation
\[\lb{f1}
 g(w)=(2\vp)g(\sqrt{1+w}),
\]
where, as above, the branch of square root satisfies $\sqrt{1+\vp}=\vp$ and it is analytic in the neighborhood of $\vp$. Substituting $w=f(z)$ into \er{f1} and applying $f^{-1}$, we obtain the Poincare equation
\[\lb{f2}
 f(2\vp z)=f(z)^2-1.
\]
{\bf Remark.} It can be shown that the functional equation $f(az)=f(z)^2-1$, $f'(0)\neq0$ admits an analytic solution if and only $a=1\pm\sqrt{5}$, see also \cite{mathforum}. 
Some general information about Poincare equations $f(az)=P(f(z))$ with polynomial or rational $P$ can be found in, e.g., \cite{F1}.

Applying Leibnitz rule to \er{f2}, we get
\[\lb{diff}
 (2\vp)^nf^{(n)}(0)=\sum_{k=0}^n\binom{n}{k}f^{(k)}(0)f^{(n-k)}(0),\ \ n\ge2,
\]
which with $f(0)=\vp$ gives the recurrent formula
\[\lb{diff1}
 f^{(n)}(0)=\frac{\sum_{k=1}^{n-1}\binom{n}{k}f^{(k)}(0)f^{(n-k)}(0)}{(2\vp)^n-2\vp},\ \ n\ge2.
\]
It is clear that $|f'(0)|=1$, $|f''(0)|=((2\vp)^2-2\vp)^{-1}\le1$. Suppose that we already proved $|f^{(k)}(0)|\le1$ for all $k=1,...,n-1$. Then, by \er{diff1}, we have
$$
 |f^{(n)}(0)|\le\frac{\sum_{k=1}^{n-1}\binom{n}{k}|f^{(k)}(0)||f^{(n-k)}(0)|}{(2\vp)^n-2\vp}\le\frac{\sum_{k=1}^{n-1}\binom{n}{k}}{(2\vp)^n-2\vp}=\frac{2^n-2}{(2\vp)^n-2\vp}\le1,\ \ n\ge2.
$$
Hence, all $|f^{(n)}(0)|\le1$ and the coefficients in Taylor expansion $|f^{(n)}(0)/n!|\le1/n!$. It leads to

{\bf Proposition.} {\it
The function $f$ can be extended to an entire function of an exponential type. Its order does not exceed $1$ (the exact order will be provided below).	
}

{\bf Remark.} Due to \er{diff1}, we conclude that all $f^{(n)}(0)>0$. Thus, $f$ is strictly increasing function for $z>0$. It is easy to check that there are negative zeros of $f$. Let $c<0$ be the first zero $f(c)=0$. Then
the values
$$
 f((2\vp)^{-n}c)=\underbrace{\sqrt{1+...+\sqrt{1}}}_n>0,\ \ n\ge1
$$
completely determine the entire function $f$. In particular, the Paris constant $2C=-c$, see Section \ref{sec0}.

{\section{Zeros of $f$ and polynomial dynamics}\lb{sec2}}

The well-known identity $\cos2z=2(\cos z)^2-1$ is similar to \er{f2}. It generates the scaled Chebyshev polynomials
$$
 \cos(2^nz)=\wt T_{2^n}(\cos z),\ \ \wt T_{2^n}(z)=\underbrace{(2z^2-1)\circ...\circ(2z^2-1)}_n=(2z^2-1)^{n\circ}.
$$
Identity \er{f2} also generates the polynomials
\[\lb{p1}
 f((2\vp)^nz)=P_n(f(z)),\ \ P_n(z)=(z^2-1)^{n\circ}.
\]
Differentiating \er{f2}, we obtain
$$
 f'(z)=\vp^{-1}f((2\vp)^{-1}z)f'((2\vp)^{-1}z)=...=\vp^{-n}f((2\vp)^{-1}z)f((2\vp)^{-2}z)...f((2\vp)^{-n}z)f'((2\vp)^{-n}z)
$$
which gives us
\[\lb{d1}
 f'(z)=\prod_{n=1}^{\iy}\vp^{-1}f((2\vp)^{-n}z),
\]
since $f'((2\vp)^{-n}z)\to f'(0)=1$. This is an analog of the Euler formula
$$
 \sin x=x\prod_{n=1}^{\iy}\cos\frac{x}{2^n}.
$$ 
For the rest, the function $f$ and $f'$ are more complex than $\cos$ and $\sin$. One of the reason is that $P_n$ can have many complex zeros. Recall that zeros of the Chebyshev polynomials are real.  Using
\[\lb{double}
 f((2\vp)^2z)=P_2(f(z))=f(z)^2(f(z)^2-2),
\]
we obtain the following

\no{\bf Proposition.}{\it If $z_0$ is a zero of $f$ of multiplicity $m$ then $(2\vp)^2z_0$ is a zero of $f$ of multiplicity $2m$.}

\begin{figure}[h]
	\begin{minipage}[h]{0.49\linewidth}
		\center{\includegraphics[width=0.9\linewidth]{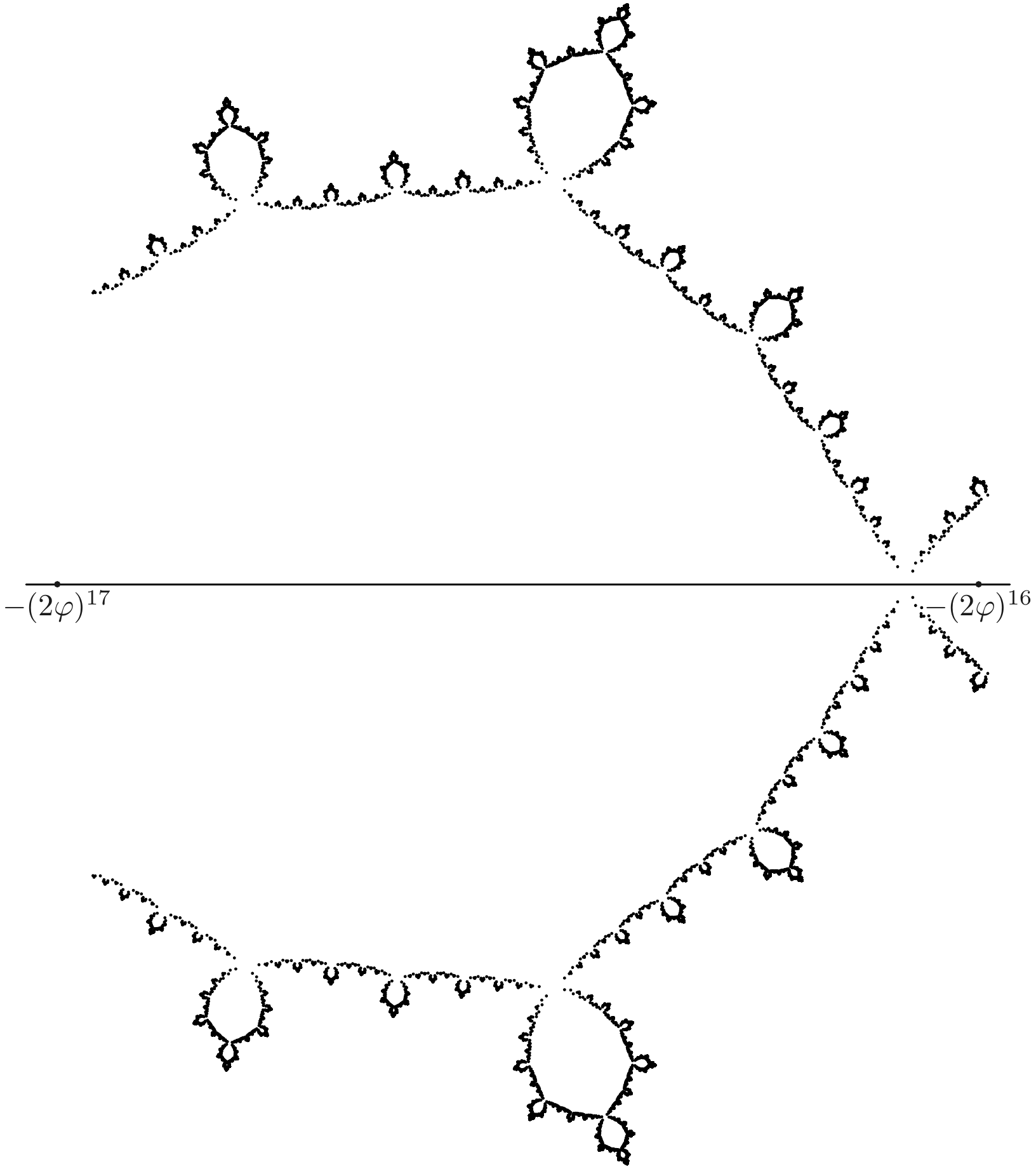} \\ (a)}
	\end{minipage}
	\hfill
	\begin{minipage}[h]{0.49\linewidth}
		\center{\includegraphics[width=0.9\linewidth]{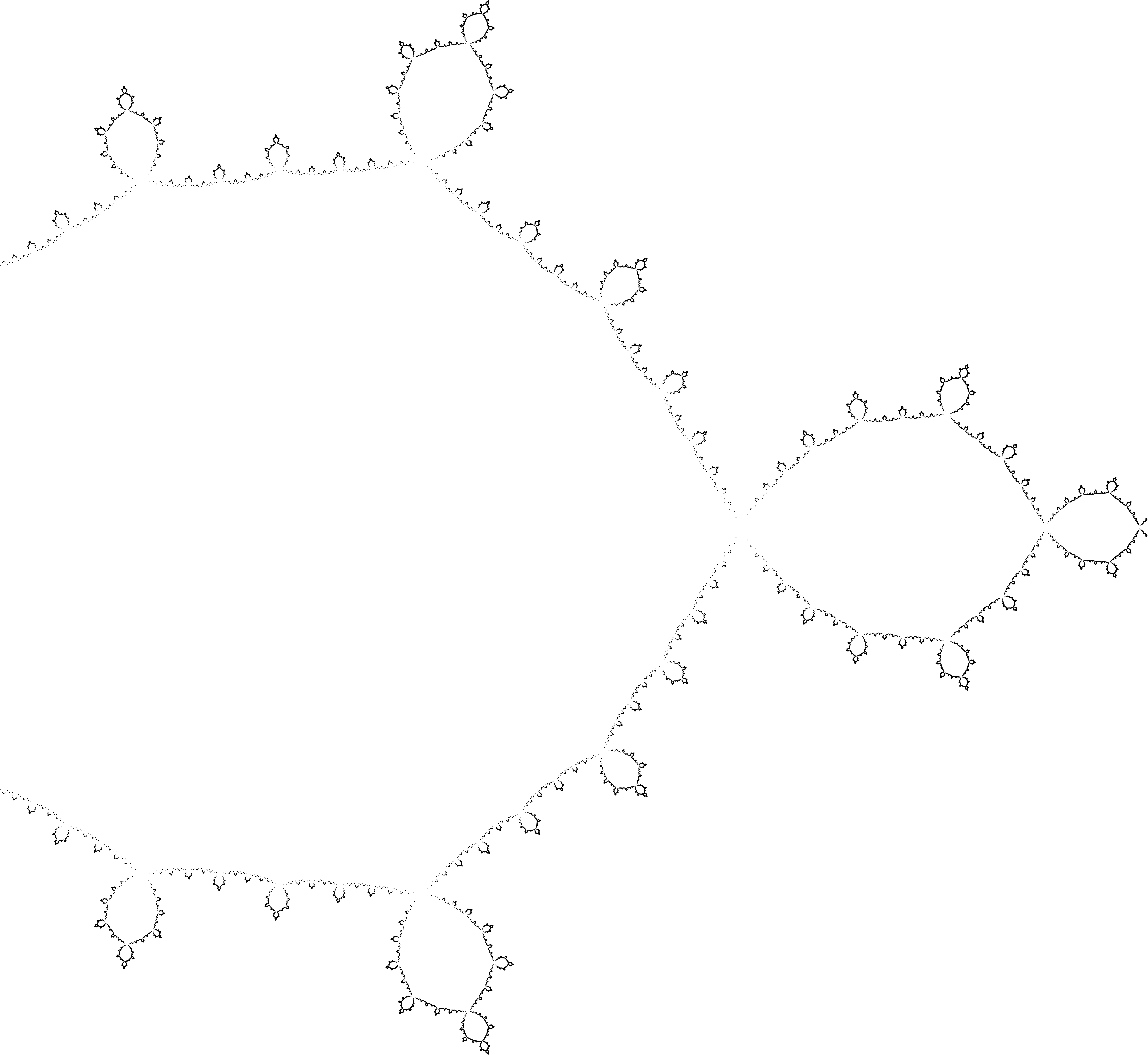} \\ (b)}
	\end{minipage}
	\caption{Approximate position of simple zeros of $f$ located in (a) $S_{16}$, (b) $S_{17}\cup S_{16}\cup S_{15}$.}
	\label{fig1}
\end{figure}

Moreover, from \er{d1} it is seen that $z_0$ is a simple zero of $f$ if and only if $f((2\vp)^{-n}z_0)\neq0$ for all $n>1$. Let $z_0$ be some simple zero of $f$. Then
\[\lb{simp1}
0=f(z_0)=P_{n}(f((2\vp)^{-n}z_0)),
\]
where $n\in\N$. It means that $f((2\vp)^{-n}z_0)$ is a {\it primitive} zero of $P_{n}$.  We call zero $y_0$ of $P_n$ primitive if $P_k(y_0)\neq0$ for all $k<n$. All primitive zeros $y$ of $P_n$, $n\ge2$ have the form
\[\lb{prim1}
 y({\bf s})=s_1\sqrt{1+s_2\sqrt{1+...+s_{n-1}\sqrt{2}}},
\]  
where ${\bf s}=(s_k)\in\{\pm1\}^{n-1}$. 
There seems to be $2^{n-1}$ distinct primitive zeros of $P_n$, $n\ge2$. 
Formula \er{prim1} can be easily implemented for numerical experiments. Using \er{simp1}, and $f(0)=\vp$, $f'(0)=1$, we conclude that
\[\lb{simp2}
 (2\vp)^{-n}z_0=y({\bf s})-\vp+O(z_0^2(2\vp)^{-2n}),
\]
for some ${\bf s}\in\{\pm1\}^{n-1}$. Denote the ring $S_n:=\{z:\ (2\vp)^n<|z|\le (2\vp)^{n+1}\}$. If $z_0$ belongs to the ring $S_{n-n_1}$ for some $n,n_1\in\N$, then 
\[\lb{simp3}
 z_0=(2\vp)^{n}(y({\bf s})-\vp)+O((2\vp)^{n-2n_1+2}).
\]
Moreover, any primitive zero $y({{\bf s}})$ of $P_n$ lying in $\vp+S_{-n_1}$, generates by \er{simp3} a simple zero $z_0$ of $f$ lying in the relative vicinity (of order $O((2\vp)^{-n_1+2})$) of $S_{n-n_1}$, when $n_1$ is large. Since such primitive zeros of $P_n$ approximate a Julia set, we obtain that large simple zeros of $f$ also create the structure similar to the Julia set. In Fig. \ref{fig1}, we draw approximate positions of simple zeroes of $f$, using primitive zeros of $P_{28}$ and $n_1=11,12,13$. 


In fact, it is possible to describe all the simple zeros of $f$ explicitly. Let $z_0$ be some simple zero of $f$. Following the arguments mentioned above, we have
$$
 f((2\vp)^{-n}z_0)=s_1\sqrt{1+s_2\sqrt{1+...+s_{n-1}\sqrt{2}}}
$$ 
for some $s_k\in\{\pm1\}$. For sufficiently large $n$, the value $(2\vp)^{-n}z_0$ lies in the small vicinity of $0$, where the inverse function $g=f^{-1}$ is well defined. Using \er{main1}, we obtain
$$
 (2\vp)^{-n}z_0=(s_1\sqrt{1+s_2\sqrt{1+...+s_{n-1}\sqrt{2}}}-\vp)\cdot\frac{2\vp}{\vp+\sqrt{1+s_1\sqrt{1+s_2\sqrt{1+...+s_{n-1}\sqrt{2}}}}}\cdot
$$
$$ 
 \frac{2\vp}{\vp+\sqrt{1+\sqrt{1+s_1\sqrt{1+s_2\sqrt{1+...+s_{n-1}\sqrt{2}}}}}}\cdot...=(s_2\sqrt{1+s_3\sqrt{1+...+s_{n-1}\sqrt{2}}}-\vp)\cdot
$$
$$
 \frac{1}{\vp+s_1\sqrt{1+s_2\sqrt{1+s_3\sqrt{1+...+s_{n-1}\sqrt{2}}}}}\cdot\frac{2\vp}{\vp+\sqrt{1+s_1\sqrt{1+s_2\sqrt{1+...+s_{n-1}\sqrt{2}}}}}\cdot...=
$$
$$
 -\frac1{\vp}\cdot\frac{1}{\vp+s_{n-1}\sqrt{2}}\cdot\frac{1}{\vp+s_{n-2}\sqrt{1+s_{n-1}\sqrt{2}}}\cdot...\cdot\frac{2\vp}{\vp+\sqrt{1+s_1\sqrt{1+s_2\sqrt{1+...+s_{n-1}\sqrt{2}}}}}....
$$
Thus
\[\lb{zero1}
 z_0((\s_n)):=z_0=-2\prod_{n=1}^{\iy}\frac{2\vp}{\vp+\s_n\sqrt{1+\s_{n-1}\sqrt{1+...+\s_1\sqrt{2}}}},
\] 
where $(\s_n)\in\S$. The set $\S$ consists of infinite sequences of $\pm1$ converging to $+1$ at infinity. Any such $z_0$ defined by \er{zero1} is simple zero of $f$. This statement can be proven using the same arguments as above. Note that the condition $(2\vp)^{-n}z_0$ is close to zero means only $(\s_{n-k}=)s_k=1$ for first $k=1,...,N$ for some $N$ depending on $z_0$ and $n$. This condition is already satisfied because the tail of sequences $(\s_n)$ consists of $+1$. Note that the only real zero $z_0((1))=-2C$, where $C$ is the Paris constant, see the end of Section \ref{sec1}.

{\section{Weierstrass-Hadamard expansion of $f$ and related formulas}\lb{sec3}}

Let $\{z_n\}$ be all simple zeros of $f$. Denote $H(z)=\prod_{n}(1-z/z_n)$. Numerical experiments show that $H$ is well defined, since the number of simple zeros located in $S_n$ grows up approximately as $A2^n$ ($A>0$), while the radius of $S_n$ is $(2\vp)^{n+1}$. The order of grows $A2^n$ is also confirmed by the exact order of $f$, see below. Recall that each simple zero $z_n$ generates zeros $(2\vp)^kz_n$ of multiplicity $2^k$, $k\ge1$.  Then, since $f$ is an entire function of the order at most $1$, it admits the Weierstrass-Hadamard expansion
$$
 f(z)=\vp e^{dz}\prod_{n=0}^{\iy}H\lt(\frac{z}{(2\vp)^{2n}}\rt)^{2^n}
$$
with some constant $d$. In fact, $d=0$, since the order of $f$ does not exceed $\r=\ln 2/\ln(2\vp)<1$. The value $\r$ can be obtained by substituting $e^{B|z|^{\r}}$ into \er{f2} which leads to $e^{(2\vp)^{\r}B|z|^{\r}}=e^{2B|z|^{\r}}$ or $(2\vp)^{\r}=2$. 
The exact value $\r$ first appeared in comments, see \cite{mathforum}. 
In principle, the orders of Poincare functions are usually easy to compute.  Thus
\[\lb{WH1}
f(z)=\vp \prod_{n=0}^{\iy}H\lt(\frac{z}{(2\vp)^{2n}}\rt)^{2^n}
\]
Substituting \er{WH1} into \er{double}, we obtain another formula for $f$ containing simple zeros only
\[\lb{primitive}
 f(z)=\sqrt{2+\vp^{-1}H((2\vp)^2z)}.
\]
Substituting \er{WH1} into \er{d1}, we obtain also the formula for derivative
\[\lb{derivative}
 f'(z)=\prod_{n=1}^{\iy}H\lt(\frac{z}{(2\vp)^{2n}}\rt)^{2^n-1}\prod_{n=1}^{\iy}H\lt(\frac{z}{(2\vp)^{2n-1}}\rt)^{2^n-1}.
\]
Finally, let us simplify \er{primitive} writing
$$
 f(2\vp z)^2-1=1+\vp^{-1}H((2\vp)^3z)
$$
or, by \er{f2},
\[\lb{primitive1}
 f(z)=1+\frac{H(2\vp z)}{\vp}.
\]
Combining \er{zero1} and \er{f2} we obtain the special closed form for $f$ formulated in
\begin{theorem}\lb{T1}
The function $f$ admits the following representation
\[\lb{explicit}
 f(z)=1+\vp^{-1}\prod_{(\s_n)\in\S}\lt(1+z\vp\prod_{n=1}^{\iy}\frac{\vp+\s_n\sqrt{1+\s_{n-1}\sqrt{1+...+\s_1\sqrt{2}}}}{2\vp}\rt).
\]	
\end{theorem}
I put Theorem here but I still do not know if this result new or not, and if there are mistakes here or not. I will be grateful if someone will indicate the original source for similar results. 

Using \er{explicit}, it is possible to obtain various identities for momenta of inverse simple zeros of $f$, e.g., the first momentum is
\[\lb{moment1}
 1=f'(0)=\sum_{(\s_n)\in\S}\prod_{n=1}^{\iy}\frac{\vp+\s_n\sqrt{1+\s_{n-1}\sqrt{1+...+\s_1\sqrt{2}}}}{2\vp}.
\]

{\bf Remark.} Direct use of \er{WH1} along with \er{zero1} leads to another special closed form
\[\lb{explicit1}
 f(z)=\vp\prod_{(\s_n)\in\S}\lt(1+\frac{z}{\vp}\prod_{n=1}^{\iy}\frac{\vp+\s_n\sqrt{1+\s_{n-1}\sqrt{1+...+\s_1\sqrt{1}}}}{2\vp}\rt),
\]
where all the multiplicities of zeros are taken into account because, e.g., $1+\s_1\sqrt{1}$ can be $0$. The analog of the first momentum formula \er{moment1} is
\[\lb{moment11}
1=f'(0)=\sum_{(\s_n)\in\S}\prod_{n=1}^{\iy}\frac{\vp+\s_n\sqrt{1+\s_{n-1}\sqrt{1+...+\s_1\sqrt{1}}}}{2\vp}.
\]

\end{document}